\documentclass[11pt,letterpaper,twoside,reqno,nosumlimits]{amsart}
\usepackage{subeqnarray,url}

\newtheorem{definition}{Definition}[section]
\newtheorem{proposition}{Proposition}[section]

\newtheorem{theorem}{Theorem}[section]
\usepackage{etoolbox}
\usepackage{amssymb}
\patchcmd{\section}{\scshape}{\bfseries}{}{}
\makeatletter
\renewcommand{\@secnumfont}{\bfseries}
\makeatother
\usepackage[toc,page]{appendix}

\usepackage[dvipsnames]{xcolor}
\patchcmd{\section}{\normalfont}{\normalfont\color{MidnightBlue}}{}{}
\patchcmd{\subsection}{\normalfont}{\normalfont\color{MidnightBlue}}{}{}

\makeatletter
\def\subsubsection{\@startsection{subsubsection}{3}%
\z@{.5\linespacing\@plus.7\linespacing}{-.5em}%
{\normalfont\bfseries}}
\makeatother

\usepackage{circuitikz}

\usepackage{fancyhdr}
\usepackage{amsmath,amsfonts,amsbsy,amsgen,amscd,mathrsfs,amssymb,amsthm,mathtools,tensor}

\usepackage[a4paper,margin=2.0cm]{geometry}
\usepackage{tikz}
\usepackage{overpic}

\usepackage{stmaryrd}

\usepackage{amsopn}

\usepackage{array} 
\usepackage{booktabs}
\usepackage{makecell}
\usepackage{caption}
\usepackage{graphicx}
\usepackage{color}
\usepackage{url}
\usepackage{epstopdf}
\usepackage{pdfsync}
\usepackage[colorlinks]{hyperref}
\usepackage{comment}
\usepackage{mathabx}
\usepackage{upgreek}
\usepackage{subfig}
\usepackage{parskip}

\usepackage{mathrsfs}
\usepackage{cleveref}
\usepackage{yfonts}

\hypersetup{
    linkcolor=blue,
}

\newlength{\fixboxwidth}
\usepackage{epsfig,amsbsy,graphicx,multirow}

\usepackage{ algorithm, algorithmic}
\renewcommand{\algorithmiccomment}[1]{\bgroup\hfill//~#1\egroup}

\usepackage[all,cmtip]{xy}

\usepackage{accents}

 \numberwithin{equation}{section}

\def\[{\begin{equation}}
\def\]{\end{equation}}

\def\RR{\mathbb{R}}

\def\N{\mathcal{N}}

\def\restrict#1{\raise-.5ex\hbox{\ensuremath|}_{#1}}

\def\<{\big\langle}
\def\>{\big\rangle}

\def\Hc{\mathcal{H}}

\def\RR{\mathbb{R}}

\definecolor{red}{rgb}{0.9, 0, 0}

\definecolor{green}{rgb}{0.0, 1.0, 0.0}

\makeatletter
\newcommand{\oset}[3][0ex]{%
  \mathrel{\mathop{#3}\limits^{
    \vbox to#1{\kern-2\ex@
    \hbox{$\scriptstyle#2$}\vss}}}}
\makeatother

\makeatletter
\newcommand{\uset}[3][0ex]{%
  \mathrel{\mathop{#3}\limits_{
    \vbox to#1{\kern-2\ex@
    \hbox{$\scriptstyle#2$}\vss}}}}
\makeatother

\usepackage{aurical}
\usepackage[T1]{fontenc}

\usepackage{tikz}
\usetikzlibrary{shadows}
\usetikzlibrary{arrows}
\usetikzlibrary{shapes}

\usepackage[numbers]{natbib}

\usepackage[foot]{amsaddr}

\begin{document}

\title[Gaussian Processes simplify differential equations]{Gaussian Processes simplify differential equations}
\author{Jonghyeon Lee$^1$, Boumediene Hamzi$^{2,5}$, Yannis Kevrekidis$^3$, Houman Owhadi$^4$}
\address{$^1$ Department of Computing and Mathematical Sciences, Caltech, CA, USA.}

\email{jlee9@caltech.edu}
\address{$^2$ Department of Computing and Mathematical Sciences, Caltech, CA, USA. }
\email{boumediene.hamzi@gmail.com}
\address{$^3$ Department of Applied Mathematics and Statistics, Johns Hopkins University, MD, USA.}
\email{yannisk@jhu.edu}
\address{$^4$ Department of Computing and Mathematical Sciences, Caltech, CA, USA. }
\email{owhadi@caltech.edu}
\address{$^5$ The Alan Turing Institute, London, UK. }

          \begin{abstract}

        In this paper we use Gaussian processes (kernel methods) to learn mappings between trajectories of distinct differential equations. Our goal is to simplify both the representation and the solution of these equations.
        We begin by examining the Cole-Hopf transformation, a classical result that converts the nonlinear, viscous Burgers' equation into the linear heat equation. We demonstrate that this transformation can be effectively learned using Gaussian process regression, either from single or from multiple initial conditions of the Burgers equation. We then extend our methodology to discover mappings between initial conditions of a nonlinear partial differential equation (PDE) and a linear PDE, where the exact form of the linear PDE remains unknown and is  inferred through Computational Graph Completion (CGC), a generalization of Gaussian Process Regression from approximating single input/output functions to approximating multiple input/output functions that interact within a computational graph. Further, we employ CGC to identify a local transformation from the nonlinear ordinary differential equation (ODE) of the Brusselator to its Poincar\'{e} normal form, capturing the dynamics around a Hopf bifurcation. We conclude by addressing the broader question of whether systematic transformations between nonlinear and linear PDEs can generally exist, suggesting avenues for future research.

\end{abstract}

\maketitle
 \section{Introduction}

\subsection{Motivations}

The study of the behavior of solutions of ordinary and partial differential equations often benefits from deciding on a convenient choice
of coordinates. This choice of coordinates may be used to “simplify” the functional expressions that
appear in the vector field in order that the essential features of the flow of an ODE, or a PDE, near a critical
point become more evident.

In the case of the analysis of an ordinary differential equation in the
neighborhood of an equilibrium point, this naturally leads to the consideration of the possibility
to remove, through a transformation, the maximum number of terms in the Taylor expansion of the vector field up to a
given order. This idea was introduced by Henri Poincar\'e in \cite{Poincare1885} and the “simplified” system is called
normal form. There have been several applications of the method of normal forms particularly in
the context of bifurcation theory where one combines the method of normal forms and
the center manifold theorem in order to classify bifurcations \cite{kuznetsov}.

For PDEs, transformations such as the Lie-Bäcklund transformation \cite{lie_backlund}, the inverse scattering transform 
 \cite{invscattering1,invscattering2} and the Miura transformation \cite{miura}  play crucial roles in simplifying complex equations and deepening our understanding of their underlying mathematical structures and their solutions. 

 In particular, transforming nonlinear PDEs to linear PDEs has generated considerable interest in research because linear PDEs are significantly easier to solve and analyze than their nonlinear counterparts. The most famous example, which we study in detail, is the Cole-Hopf transformation \cite{cole,hopf} because it describes the global dynamics of a nonlinear PDE, the viscous Burgers equation, in terms of a linear PDE, the heat equation. The Cole-Hopf transformation has found applications in fields as varied as quantum mechanics \cite{scarfone2002} and stochastic analysis \cite{cirant2015}. 
 
 

 Previous work on linearizing nonlinear dynamics are divided into two broad categories: 'formula-based' (analytic) approaches and data-driven approaches, to which our paper belongs. The former category includes direct derivations in the style of Cole-Hopf \cite{wuchen2021}. The most common approach in the latter category involves Koopman operators, often combined with neural networks, which involves mapping a nonlinear system to an infinite dimensional linear system whose output is encoded by a Koopman matrix (a finite dimensional approximation of the infinite dimensional space) \cite{Mezic05, bookmezic}. However, Gaussian processes (GPs) have the advantage of being simple to implement and offer superior convergence results. GPs and reproducing kernel Hilbert spaces (RKHS) also provide mathematical foundations for analyzing dynamical systems \cite{ bhcm11,bhcm1,lyap_bh,bhks,hamzi2019kernel, bh2020b,klus2020data,ALEXANDER2020132520,bh12,bh17,hb17, mmd_kernels_bh, akian2022learning, HAMZI2023128583, 5706920,   yk1, bh_sparse_kfs, BHPhysicaD}, surrogate modeling \cite{santinhaasdonk19}, and even analyzing neural networks \cite{owhadi20,smirnov2022mean}. Most importantly to our paper, Gaussian processes are especially powerful because they condense the problem of linearizing nonlinear PDEs, when the linear PDE we map our nonlinear PDE to is known, to a simple regression formula.

Our work also explores the application of transformations to find local transformations between nonlinear dynamical systems and their normal forms. We consider a simple case in which we examine a trajectory generated by the Brusselator and aim to map it to the corresponding trajectory of its Poincar\'e normal form that exhibits a Hopf bifurcation. Such transformations exist due to topological equivalence and can be approximated by Taylor series approximations \cite{kuznetsov} but kernels allow us to learn them using only the data points taken from the trajectory of the dynamical system without assuming knowledge of the underlying equations that generate them.
 
\par
\par

\subsection{Outline of the paper}

The structure of the paper can be summarized as follows. In Section 2, we show that the Cole-Hopf transformation between a pair of initial conditions of the Burgers and heat equations can be learned by kernel regression. We also illustrate the impact the choice of kernel can have on the accuracy and describe an effective method for learning the kernel. We then illustrate that Cole-Hopf can also be learned using multiple initial conditions of the Burgers' equation.
In Section 3, we apply our method to another pair of solutions emanating from given initial conditions of nonlinear and linear PDEs, before using a nonlinear optimization problem to learn the transformation even when the exact form of the linear PDE is not known. Section 4 uses a similar  framework to learn a map between a trajectory of the Brusselator and the corresponding trajectory of the normal form equations where all the parameters of the dynamical systems are known. Finally, we outline a process for learning further transformations between nonlinear and linear PDEs via Gaussian process regression.

\section{Cole-Hopf transformation}

\subsection{Gaussian Process Regression for the Cole-Hopf transformation}

\begin{figure}[tbh]
\centering
\subfloat{\label{fig:test1}
\includegraphics[width=.45\textwidth]{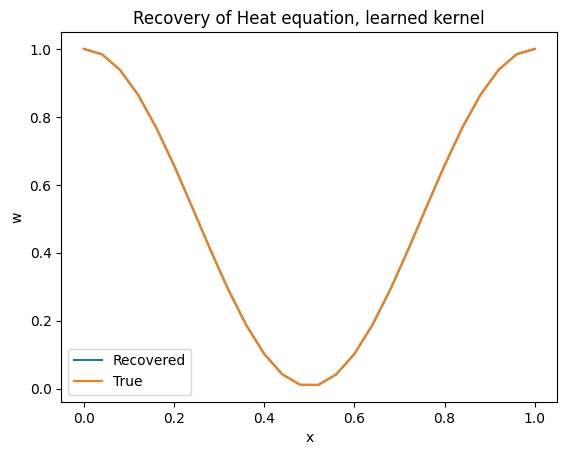}
}
\hfill
\subfloat{\label{fig:test2}
\includegraphics[width=.45\textwidth]{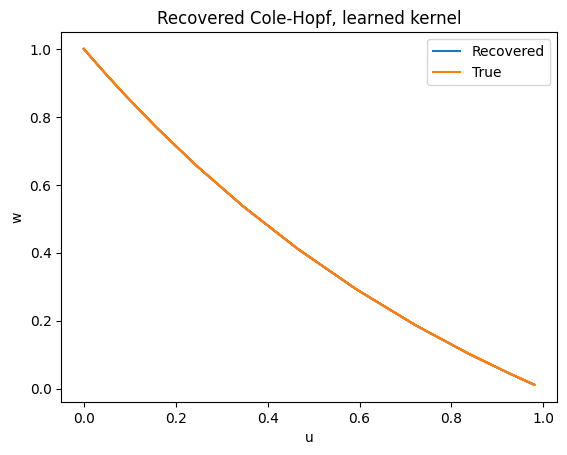}
}
\hfill
\caption{Cole-Hopf transformation: true $D^\dagger (u)=w$ and learned $D(u)$ against $x$ (left) and $u$ (right), $N=25$ points between $x=0$ and $x=1$}
\label{chplot}
\end{figure}

In this section, we explore a specific version of the Cole-Hopf transformation (see Figures \ref{chplot} and \ref{chplot2}) that connects the initial condition of the Burgers equation with that of the heat equation\footnote{See Appendix A.3 for more details about the Cole-Hopf transformation}:

\begin{equation}
    w_0 = D^\dagger (u_0) := \frac{e^{-\frac{1}{2\nu}u_0}-e^{-\frac{1}{2\nu}}}{1-e^{-\frac{1}{2\nu}}},
    \label{eq:truecolehopf}
\end{equation}
where $w_0$ is an initial condition of the heat equation, and $u_0$ is the antiderivative of the initial condition of the Burgers equation, $v_0$.

Note that the transformation maps the value of $u_0$ at a point $x$ to the value of $w_0$ at the same point in a manner independent of $x$. Thus, the transformation can be interpreted both as a function that takes $u_0$ to $w_0$, and as an operator that pointwise maps the entire initial condition $u_0(x)$ to $w_0(x)$.
In its current form, the transformation acts on the values of the  antiderivative of $v_0$  at each point $x$, which is advantageous as it circumvents the curse of dimensionality by reducing the problem to a lower-dimensional space. This method allows for high-accuracy identification of the transformation by leveraging specific structural knowledge. Alternatively, the transformation could be designed to act on the entire function $x \mapsto v_0(x)$. While this approach is more general and does not rely on detailed structural information, it suffers from the curse of dimensionality due to the high-dimensional nature of the function space. The current formulation is essentially a specialized case of this broader approach, employing a  specific operator-valued kernel that simplifies the problem to a pointwise transformation.

We start by assuming that we know the initial condition of the Burgers equation $v_0$, as well as its primitive $u_0 = Pv_0$, where $P$ denotes the integral operator $Pv_0 := \int^x_0 v_0 (s) ds$. Then the solution of the Burgers equation at time $t_1 = h$, where $h$ is small, can be approximated using a simple forward time solver based on the Euler method:

\begin{equation}
    v_1 = Bv_0 := v_0 + h (\nu \partial^2_x v_0-v_0 \partial_x v_0).
\label{burgerssolver}
\end{equation}

Similarly, the forward Euler method can also be used to evolve the initial condition of the heat equation $w_0$ by a single time step:

\begin{equation}
    w_1 = Hw_0 := w_0 + h \nu \partial^2_x w_0.
\label{heatsolver}
\end{equation}

\begin{figure}[ht]
    \centering
    \includegraphics{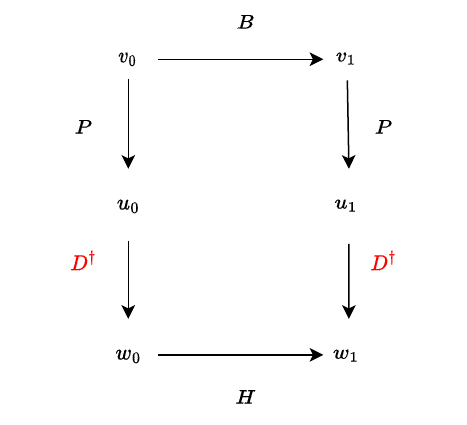}
    \caption{Relationship between the solution of the heat equation $w$ and the solution of the Burgers equation $v$; note that $w_1 = HDPv_0 = DPBv_0$}
    \label{ColeHopfdiagram}
\end{figure}

From the diagram, it is clear that $w_0 = D^\dagger u_0 = D^\dagger Pv_0$ and $w_1 = Hw_0 = D^\dagger u_1 = D^\dagger Pv_1 = D^\dagger PBv_0$. Our approach to finding our unknown function $D^\dagger$ is to replace it by a centered Gaussian process (GP) with kernel  $K$ and to learn the conditional expectation (or equivalently, the MAP estimator) of the GP given data. Thus, the problem of recovering $D^\dagger$ and $w_0$ from $v_0$ becomes that of finding the candidate function $D \in \mathcal{H}_K$ that minimizes the squared RKHS norm $||D||^2_K := [w_0 , K^{-1} w_0]$ with respect to $K$ under the following constraints, where $N$ is the number of collocation points and $\mathcal{H}_K$ is the RKHS associated with the kernel $K$:

\begin{equation}
\begin{aligned}
\min_{D \in \mathcal{H}_K} \quad & ||D||^2_K \\
\textrm{s.t.} \quad & DPv_0 (x_i) = w_0 (x_i), \; i= 1, \dots, N \\
\quad & DPBv_0 (x_i) = Hw_0 (x_i), \; i= 1, \dots, N \\ 
\quad & D(0) =1 \\
\quad & D(1) = 0 
\end{aligned}.
\label{eq:colehopfminimizer}
\end{equation}

We make a few pertinent observations about (\ref{eq:colehopfminimizer})\footnote{(\ref{eq:colehopfminimizer}) is a discretized version of a homological equation}.
First, we note that if there exists a $D$ that satisfies the first two constraints in \ref{eq:colehopfminimizer}, then $aD+b$ and $D+cu$, where $a,b,c \in \mathbb{R}$ are constants, also satisfy those constraints due to the second derivatives in (\ref{burgerssolver}) and (\ref{heatsolver}). Hence we require two uniqueness constraints to ensure that $D$ has a unique, non-zero solution. We have set them 
as $D(0)=1$ and $D(1)=0$ in (\ref{eq:colehopfminimizer}); this means that we are searching for a version of the Cole-Hopf transformation which maps the constant function $0(x)=0 \forall x \in \mathbb{R}$ to the constant function $1(x)=1 \forall x \in \mathbb{R}$, and maps the constant function $1(x)$ to the constant function $0(x)$. Because of the superposition property of the heat equation (any linear combination of solutions is also a solution of the heat equation), we may change the right-hand side of those two constraints to any two other constant functions and then the solution of (\ref{eq:colehopfminimizer}) would be another solution of the heat equation satisfying different initial conditions. Second, if no function between $u_0$ and $w_0$ were to exist that satisfied the constraints, then the recovered 'map' would become increasingly irregular and the empirical RKHS norm of $D$ would go to infinity as the number of collocation points increases. Conversely, if there exists a function $D(u_0)=w_0$, the computed RKHS norm of the approximation would still increase, but it would eventually converge to the RKHS norm of the true function as $N$ increases. Thirdly, since the Cole-Hopf transformation is global, we can replace $v_0$ and $w_0$ with $v_{t_j}$ and $w_{t_j}$, the solution of the Burgers and heat equations at an arbitrary time $t_j$, and the problem (\ref{eq:colehopfminimizer}) would still be valid.

Equivalently, we may also define a collection of functionals $\phi_i : F \rightarrow HFPv_0(x_i) - FPBv_0 (x_i)$ with the property that $\phi_i (D^\dagger)=0$ and reformulate (\ref{eq:colehopfminimizer}) as the following problem, where $u_0 (x_i) = \int^{x_i}_0 v_0 (x) dx$ are our new collocation points:

\begin{equation}
\begin{aligned}
\min_{D \in \mathcal{H}_K} \quad & ||D||^2_K \\
\textrm{s.t.} \quad & \phi_i (D) = 0, \; i =1, \dots, N\\
\quad & D(0) = 1 \\
\quad & D(1) = 0 
\end{aligned}.
\label{eq:colehopfminimizer2}
\end{equation} 

We may expand the constraints as follows, where $\delta_{u_0 (x_i)} (u)$ is the Dirac delta distribution centered at $u_0 (x_i)$\footnote{$\int_{\mathbb{R}} \delta_{u_0 (x_i)} f(u) du = f(u_0 (x_i))$ for a function $f$}:

 \begin{equation}
\begin{aligned}
\quad & \phi_i (D) =  0 = HDPv_0 (x_i) - DPBv_0 (x_i) \\
\quad & = \int_{\mathbb{R}} \delta_{u_0 (x_i)} D(u)  du + h \nu \int_{\mathbb{R}} \delta_{u_0 (x_i)} D(u) \cdot \frac{\partial^2}{\partial x^2} du - \int_{\mathbb{R}} \delta_{u_1 (x_i)} D(u) du \\
 \quad & D(1) = \int_{\mathbb{R}} \delta_1 (u) D(u) du = 0\\
\quad & D(0) = \int_{\mathbb{R}} \delta_0 (u) D(u) du= 1 \\
\end{aligned}.
\label{eq:chsimplified2}
\end{equation}

These expansions lead us to define the following distributions:
\begin{equation}
\begin{aligned}
\quad & \tilde{\phi}_{1} (u) = \delta_0 (u)\\
\quad & \tilde{\phi}_{i+1} (u)= \delta_{u_0 (x_i)} (u)+ h \nu \delta_{u_0 (x_i)} (u) \cdot \frac{\partial^2}{\partial x^2}-\delta_{u_1 (x_i)} (u)  , \; i = 1,\dots ,N  \\
\quad & \tilde{\phi}_{N+2} (u)= \delta_1 (u) \\
\end{aligned}.
\label{eq:chsimplfied}
\end{equation}

The constraints in (\ref{eq:colehopfminimizer2}) may therefore be combined to give the following regression problem which is linear in $D$:

\begin{equation}
\begin{aligned}
\min_{D \in \mathcal{H}_K} \quad & ||D||^2_K \\
\textrm{s.t.} \quad & \tilde{\phi}(D) = Y \\
\end{aligned},
\label{eq:colehopfminimizer3}
\end{equation}

where $Y =  \begin{pmatrix}
           1 \\
           0 \\
           \vdots \\\
           0
         \end{pmatrix} \in \mathbb{R}^{N+2}$ and $\tilde{\phi}(D)_i = [\tilde{\phi}_i,D] := \int_{\mathbb{R}} \tilde{\phi}_i (u) D(u) du $.
In other words, each entry of $Y$ is the action of a bounded linear functional $\tilde{\phi}_i$ on the unknown function $D$, where $\tilde{\phi}_i $ depends on $u_0 (x_i)$. Relaxing the constraints in  \eqref{eq:colehopfminimizer3} leads to 
\begin{equation}
\min_{D \in \mathcal{H}_K} \quad  ||D||^2_K +\frac{1}{\lambda}\big| \tilde{\phi}(D) - Y \big|^2,
\label{eq:colehopfminimizer3relaxed}
\end{equation}
whose minimizer is given by the generalized representer theorem
\begin{equation}
    D(u) = K(u,\tilde{\phi})(K(\tilde{\phi},\tilde{\phi})+\lambda I)^{-1} Y.
    \label{eq:representer},
\end{equation}
In the formula \eqref{eq:representer}, $K(u,\tilde{\phi})$ is the $N$-vector with entries $K(u,\tilde{\phi}_i)$ and $K(\tilde{\phi},\tilde{\phi})$ is the $N \times N$ matrix with entries $[\tilde{\phi}_i, K \tilde{\phi}_j]$. The nugget  $\lambda> 0$ is taken to be small enough so that \eqref{eq:representer} remains a good approximant of the minimizer of the constrained problem \eqref{eq:colehopfminimizer3}, while alleviating the possible ill-conditioning of $K(\tilde{\phi},\tilde{\phi})$.

We note that there is a certain amount of  error that results from solving the heat and Burgers equations numerically. To avoid the  error that arises from the numerical solvers of the differential equations, we take the limit as the time step $h$ approaches 0 in the formula $DPBv_0 = HPBv_0$ and derive the following differential equation:

\begin{equation}
    \frac{1}{2} \frac{\text{d}^2}{\text{d}u^2}D(u) + \nu \frac{\text{d}}{\text{d}u} D(u) =0.
\end{equation}

Hence our new regression problem becomes

\begin{equation}
\begin{aligned}
\min_{D \in \mathcal{H}_K} \quad & ||D||^2_K \\
\textrm{s.t.}\quad & D(0) =1 \\
\quad & \frac{1}{2 } \frac{\text{d}^2}{\text{d} u^2} D(u_0 (x_i))+ \nu \frac{\text{d}}{\text{d} u} D(u_0 (x_i))= 0, \; i=1,\dots,N \\
\quad & D(1) = 0 
\end{aligned}.
\label{eq:colehopfdiff}
\end{equation}

This is still a regression problem in the form of (\ref{eq:colehopfminimizer2})
 which is linear in the constraints and whose solution is given by (\ref{eq:representer}), where the $\tilde{\phi}_i (u)$ are now:

\begin{equation}
\begin{aligned}
\quad & \tilde{\phi}_1 (u)= \delta_0 (u) \\
\quad & \tilde{\phi}_{i+1} (u)= \frac{1}{2} \delta_{u_0 (x_i)}(u) \cdot \frac{\text{d}^2}{\text{d} u^2} + \nu \delta_{u_0 (x_i)}(u) \cdot \frac{\text{d}}{\text{d} u} , \; i = 1,\dots ,N \\
\quad & \tilde{\phi}_{N+2} (u) = \delta_1 (u)\\
\end{aligned}.
\label{eq:colehopfdiffsimplified}
\end{equation}

 We note that (1) this is the same ODE which was involved in the analytic derivation of the Cole-Hopf transformation and (2) any well-posed linear differential equation with sufficient uniqueness conditions can be solved as a quadratic minimization problem in the mold of (\ref{eq:colehopfdiff}).

\subsection{Learning the kernel parameters}

We use a single initial condition $u_0$ of $u$ as our training data and compare two approaches for learning Cole-Hopf. In Method A\footnote{Method A is a variant of Kernel Flows; see Appendix A.2} (Figure \ref{chplot}), we learn the parameter $\theta$ of the kernel by minimizing the following loss function:

\begin{equation}
\rho(\theta) := \frac{1}{N} \sum^N_{j=1} \left(1- \frac{Y^{-j,T} (K^{-j}_\theta (\tilde{\phi},\tilde{\phi})+\lambda I)^{-1} Y^{-j}}{Y^T (K_\theta (\tilde{\phi},\tilde{\phi})+\lambda I)^{-1} Y}\right).
\end{equation}

The matrix $K^{-j}_\theta (\tilde{\phi},\tilde{\phi}) \in \mathbb{R}^{(N+1) \times (N+1)}$ is formed by removing the row and column of $K_\theta (\tilde{\phi},\tilde{\phi})$ corresponding to the $j$th data point of $u_0 (x)$, and $Y^{-j} = (1,0,\dots,0)^T \in \mathbb{R}^{N+1}$. The intuition behind the function $\rho (\theta)$ is that the kernel parameter(s) is considered good if the interpolant does not change much when we remove one point from the training data.

\begin{figure}[tbh]
\centering
\subfloat{\label{fig:test3}
\includegraphics[width=.45\textwidth]{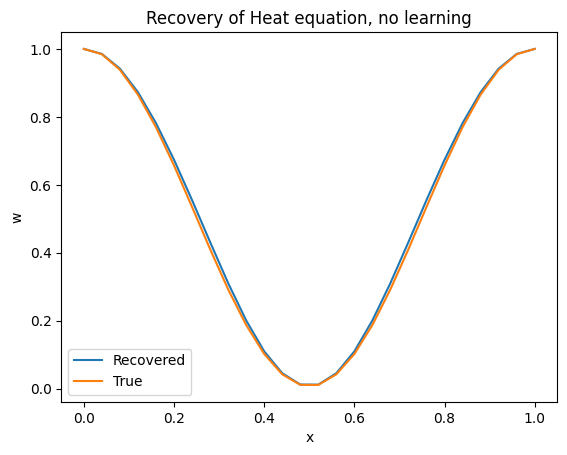}
}
\hfill
\subfloat{\label{fig:test4}
\includegraphics[width=.45\textwidth]{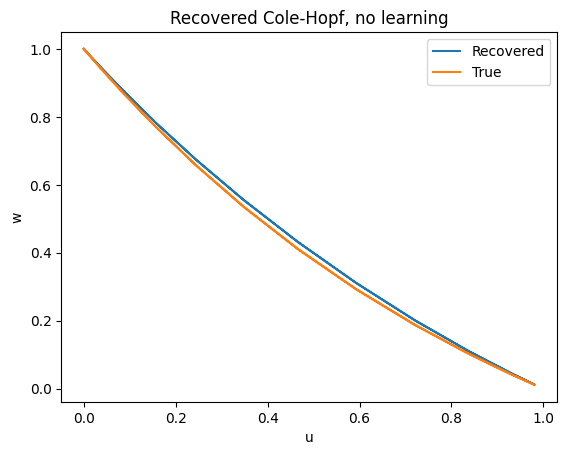}
}
\hfill
 \caption{Cole-Hopf transformation: true $D^\dagger (u)=w$ and learned $D (u)$ against $x$ (left) and $u$ (right), $N=25$ points, kernel parameter not learned}
\label{chplot2}
\end{figure}

  In Method B (Figure \ref{chplot2}), we do not learn the kernel and use the parameter $\theta = 1$. We report the relative operator loss between the true Cole-Hopf transform $D^\dagger$ and its Gaussian process approximation $D$:

\begin{equation}
    \mathcal{L} (D,D^\dagger) = \frac{||D-D^\dagger||_2}{||D^\dagger||_2}.
\label{eq:loss}
\end{equation}

In each of our examples, we use the single-parameter Matérn-2.5 kernel because it is a universal kernel that does not require many assumptions about the structure of the data:

\begin{equation}
    K_\theta (x,y) = \left(1+ \sqrt{5} \frac{|x-y|}{\theta} + \frac{5}{3} \frac{|x-y|^2}{\theta^2}\right) \exp \left(\frac{-\sqrt{5}|x-y|}{\theta}\right).
\end{equation}

 We have displayed the $\mathcal{L}(D,D^\dagger)$ in Table \ref{colehopftable} with various values of $N$ under Methods A and B.

\begin{table}[ht]
\begin{tabular}{lllllll}
\hline
\multicolumn{1}{c}{$N$} & \multicolumn{1}{c}{25}        & \multicolumn{1}{c}{50} & \multicolumn{1}{c}{100} & \multicolumn{1}{c}{200} & \multicolumn{1}{c}{400} & \multicolumn{1}{c}{800}       \\ \hline
(A) Learning            & 2.9675e-4 & 7.6794e-5              & 1.9450e-5               & 4.8647e-6               & 1.1890e-6               & 2.6892e-7 \\
(B) No learning         & 1.9232e-2 & 5.2601e-3              & 1.3532e-3               & 3.4103e-4               & 8.4719e-5               & 2.0338e-5
\end{tabular}
\caption{Relative $L^2$ error for Cole-Hopf with and without learning the kernel with $N$ points}
\label{colehopftable}
\end{table}

Table \ref{colehopftable} and Figures \ref{chplot} show that when we learn the kernel, Method A outperforms Method B in mapping the initial condition of the Burgers equation $v_0 (x) = \frac{28 \nu \pi \sin (\pi x)}{7+3.2 \cos (\pi x)}, \nu = 0.5$ to the initial condition of the corresponding heat equation $w_0 (x)$. However, as the number of data points increases, the Gaussian process approximation $D$ converges rapidly to the true function $D^\dagger$ in both cases.

\subsection{Results with multiple initial conditions}

We proceed to learn the Cole-Hopf transformation with four initial conditions of the Burgers' equation. The procedure is identical to that in the case with one initial condition, and we show that the transformation can be learned for a wider range of values of $v_0 (x)$ than in Section 2.1. First, we take 101 data points each from the following functions
\begin{equation}
\begin{aligned}
\quad & v_{0,1} (x) = 5+3x,  \ -2.5 \leq x \leq -1.5 \\
\quad & v_{0,2} (x) =  5\cos(x) +2, \; 0 \leq x \leq 1 \\
\quad & v_{0,3} (x) = \frac{\exp(\frac{x}{3})}{100}, \; 15 \leq x \leq 16 \\
\quad & v_{0,4} (x) = x, \; 10 \leq x \leq 11\\
\end{aligned}.
\end{equation}

\begin{figure}[tbh]
\centering
\includegraphics[width=.45\textwidth]
{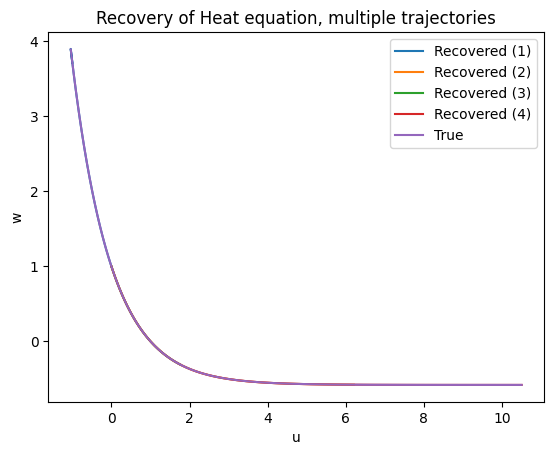}
\hfill
\caption{Cole-Hopf transformation for 4 trajectories in the same plot}
\label{fig:multipletrajectories}
\end{figure}

We use the notation $u_{0,k} := \int^x_0 v_{0,k} (s) ds, k = 1, \dots, 4$. Then solving the regression problem 

\begin{equation}
\begin{aligned}
\min_{D \in \mathcal{H}_K} \quad & ||D||^2_K \\
\textrm{s.t.}\quad & D(0) =1 \\
\quad & \frac{1}{2} \frac{\text{d}^2}{\text{d} u^2} D(u_{0,k} (x_i))+ \nu \frac{\text{d}}{\text{d} u} D(u_{0,k} (x_i))= 0, \; i=1,\dots,N, \; k = 1,\dots,4 \\
\quad & D(1) = 0 
\end{aligned}.
\label{eq:colehopfmultiple}
\end{equation}

recovers the Cole-Hopf transformation for each trajectory. Our results, shown in \ref{fig:multipletrajectories}, confirm that the Cole-Hopf transformation can be recovered with a high degree of accuracy in a larger interval in space than observed in Section 2.2.

\section{Linearizing a nonlinear PDE: a second example}

\begin{figure}[tbh]
\centering
\subfloat{\label{fig:test5}
\includegraphics[width=.45\textwidth]{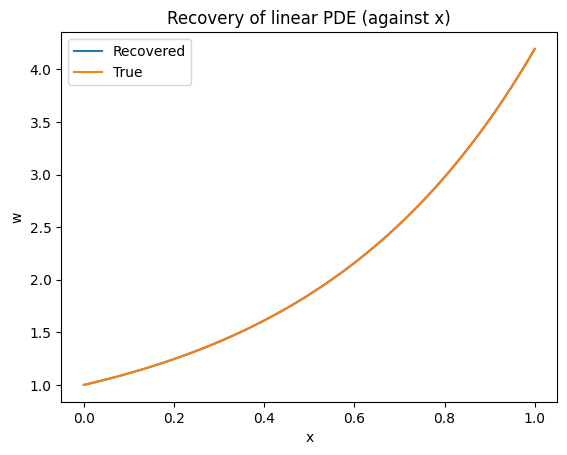}
}
\hfill
\subfloat{\label{fig:test6}
\includegraphics[width=.45\textwidth]{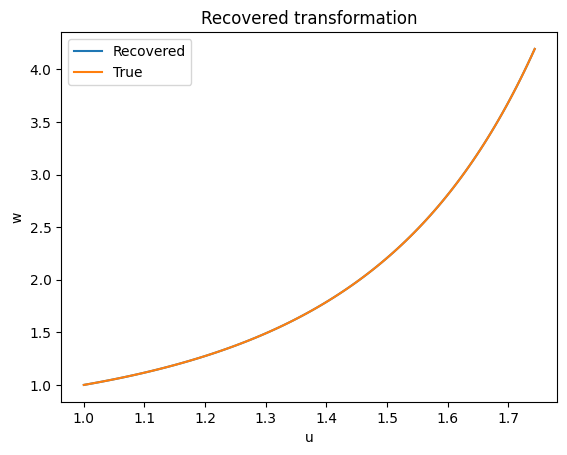}
}
\hfill
\caption{True $G^\dagger(u)=w$ and learned $G(u)$ against $x$ (left) and $u$ (right) with $N=100$ points between $x=0$ and $x=1$}
\label{linplot}
\end{figure}

We examine another nonlinear PDE that can be transformed directly into a linear PDE. The transformation, which we learn to a high level of accuracy, is depicted in Figure \ref{linplot}.

Consider the first-order nonlinear PDE

\begin{equation}
    u_t = u_x - \frac{1}{u^2}.
\label{eq:nonlinear}
\end{equation}

Suppose we wish to map an initial condition of $u$ to an initial condition of the linear PDE

\begin{equation}
    w_t = w_x - w.
\label{eq:linear}
\end{equation}

Then the Gaussian process regression problem for the transformation may be derived in an identical fashion to the previous section. We write $N$ and $M$ to symbolize forward solvers for (\ref{eq:nonlinear}) and (\ref{eq:linear}):

\begin{equation}
    u_1 = Nu_0 := u_0 + h \left( \partial_x u_0-\frac{1}{u^2_0}\right) 
\end{equation}

\begin{equation}
    w_1 = Mw_0 := w_0 + h (\partial_x w_0-w_0) 
\end{equation}

\begin{figure}[ht]
    \centering
    \includegraphics{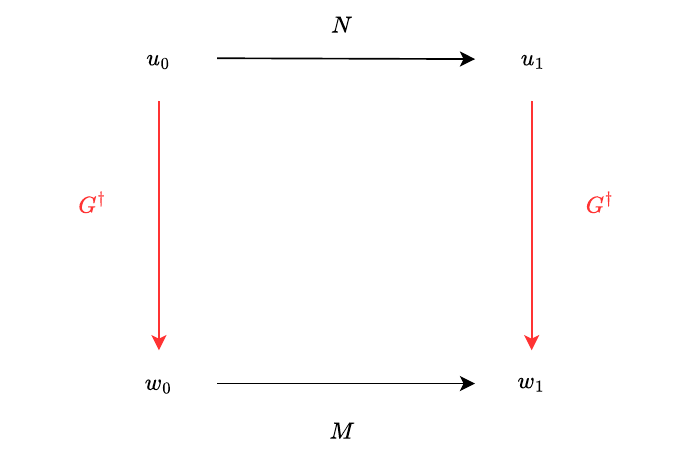}
    \caption{Relationship between the linear equation $w$ and nonlinear equation $u$}
    \label{seconddiagram}
\end{figure}

Diagram \ref{seconddiagram} visualizes the relationship between $u$ and $w$. As with the Cole-Hopf example, we can write $w_1$ in two different ways: $w_1 = MG^\dagger u_0$ and $w_1 = G^\dagger N u_0$. This time, we directly take the limit as $h \rightarrow 0$ and express the candidate function $G$ in terms of a linear differential equation: 

\begin{equation}
\frac{1}{u^2} \frac{\text{d}}{\text{d}u}G(u) = G (u)
\label{eq:lindiffeq}
\end{equation}

Once again, the solution of (\ref{eq:lindiffeq}) can be expressed as the solution of a quadratic kernel optimization problem with solution (\ref{eq:representer}), although we only need one uniqueness constraint (which we take to be $G (1) = 1$) as our PDE is first-order in space:

\begin{equation}
\begin{aligned}
\min_{G \in \mathcal{H}_K} \quad & ||G||^2_K \\
\textrm{s.t.}\quad & G(1) =1 \\
\quad & \frac{1}{u_0 (x_i)^2} \frac{\text{d}}{\text{d} u} G (u_0 (x_i)) -G(u_0 (x_i))= 0, i=1,\dots,N\\
\end{aligned}
\label{eq:secondopt}
\end{equation}

This is still a regression problem in the form of (\ref{eq:colehopfminimizer2})
 which is linear in $G$ and whose relaxed solution is given by

 \begin{equation}
    G(u) = K(u,\tilde{\phi})(K(\tilde{\phi},\tilde{\phi})+\lambda I)^{-1} Y,
    \label{eq:representer2}
\end{equation}
 
 where $Y =  \begin{pmatrix}
           1 \\
           0 \\
           \vdots \\\
           0
         \end{pmatrix} \in \mathbb{R}^{N+1}$ and the $\tilde{\phi}_i (u)$ are now:

\begin{equation}
\begin{aligned}
\quad & \tilde{\phi}_1 (u)= \delta_1 (u) \\
\quad & \tilde{\phi}_{i+1} (u)= \frac{1}{u_0 (x_i)^2} \delta_{u_0 (x_i)}(u) \cdot \frac{\text{d}}{\text{d} u} +  \delta_{u_0 (x_i)}(u)  , \; i = 1,\dots ,N \\
\end{aligned}
\label{eq:seconddiffsimplified}
\end{equation}

In Figure \ref{linplot}, we have used the representer formula with the Matérn-2.5 kernel to learn the map between $u_0 (x) = (3\log ( \cosh (-x) \exp(x))+1)^{\frac{1}{3}}$ and $w_0 (x)$. The high accuracy of the learned map, even though we have not learned the parameter of the kernel in this case, means that it is visually indistinguishable from the true transformation 

\begin{equation}
G^\dagger (u) = e^{\frac{u^3-1}{3}}
\end{equation} in the plots.

\subsection{Learning the linear PDE}

\begin{figure}[tbh]
\centering
\subfloat{\label{fig:test7}
\includegraphics[width=.45\textwidth]{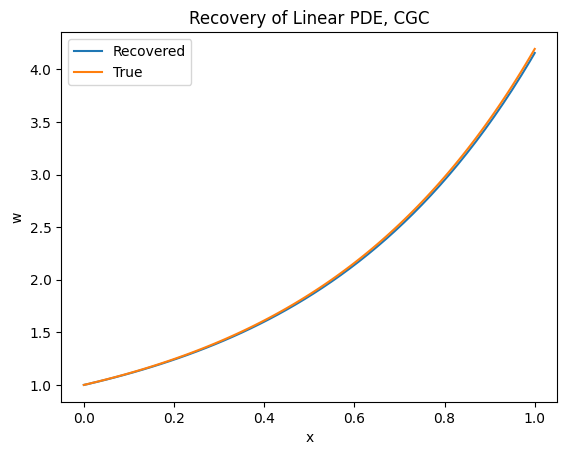}
}
\hfill
\subfloat{\label{fig:test8}
\includegraphics[width=.45\textwidth]{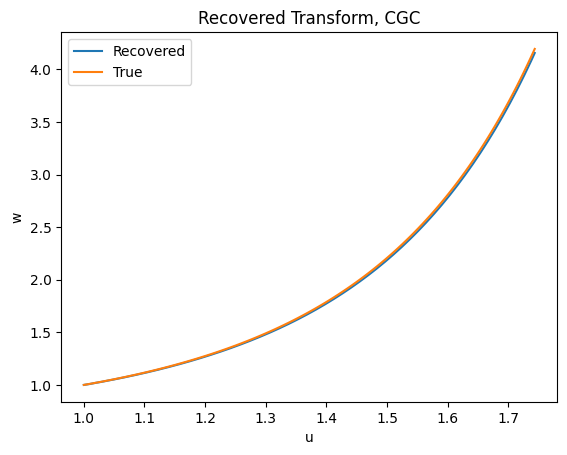}
}
\hfill
\caption{True $G^\dagger (u)=w_0$ and learned $G(u)$ against $x$ (left) and $u_0$ (right), where $a \approx -0.9974$ has been learned by CGC}
\label{cgcplot}
\end{figure}

So far, we have assumed that we know the exact equation of the linear PDE which we map our nonlinear PDE to. We proceed to explore what happens if instead, we only know that the equation of the linear PDE takes the form

\begin{equation}
    w_t = w + a^\dagger w_x ,
\end{equation}

where $a^\dagger \in \mathbb{R}$ is an unknown real number. In that case\footnote{This case can also be viewed as a way of introducing uncertainty in the target PDE since the linear PDE is no longer exactly $w_t=w-w_x$ when we initialize the algorithm for minmizing the loss function.}, we replace $a^\dagger$ with a Gaussian random variable $a$ with distribution $\mathcal{N} (0,\Gamma)$\footnote{We can generalize this to the case where the coefficient is a function $a(x)$, although we would use a non-constant kernel in that situation.}, hence the squared RKHS norm $||a||^2_\Gamma$ is simply $\frac{a^2}{\Gamma^2}$; here we take $\Gamma = 1$ (that is, $\Gamma$ is a constant kernel). Then we find $G$ and $a$ via a nonlinear optimization problem in which we obtain the best Gaussian process approximation of our unknown functions by finding the minimizer of the following loss:

\begin{equation}
    \mathcal{L}(f,Z) =  ||G||^2_K + ||a||^2_\Gamma  + \lambda_1 \mathcal{L}_1 (f,Z) + \lambda_2 \mathcal{L}_2 (Z) + \lambda_3 |G(1)-1|^2 ,
\label{eq:secondproblem}
\end{equation}

where $f =  (G,a)$ and $Z = (Z_{\cdot, 0},Z_{\cdot,1},Z_{\cdot,2})$ is the $N \times 3$ data array containing our approximations of $(u,G,\frac{\text{d}G}{\text{d}u})$\footnote{In the case of $u$, $Z_{\cdot,0}$ would be array containing the true value of $u$ at space points $x_j \in [0,1]$} at a number of collocation points. In the expressions below, we use $| \cdot |$ to denote the Euclidean norm of the arrays and $/$ to symbolize elementwise division of vectors.


\begin{equation}
  \mathcal{L}_1 (f,Z) =   |G(Z_{\cdot,0})-Z_{\cdot,1}|^2
\end{equation}

is a loss which forces $Z_{\cdot,1}$, our array of data points for the approximation of $w$, to be the output of $G(Z_0)$.


\begin{equation}
    \mathcal{L}_2 (Z) = |Z_{\cdot,1} + a Z_{\cdot,2}/Z_{\cdot,0}^2|
\end{equation}

 is a loss which enforces the differential equation (\ref{eq:lindiffeq}).

 Here, $\lambda_1$ is $\frac{1}{\lambda}$, where $\lambda$ is the regularization constant from the previous sections; $\lambda_2$ and $\lambda_3$ are typically learned by cross-validation, although in our case, we have determined them such that no term in $\mathcal{L}(f,Z)$ becomes too dominant.

Figure \ref{cgcplot} illustrates the recovered transformation $G(u)$ where we have learned $a$ by minimizing (\ref{eq:secondproblem}); our initial guesses for $a$ and $G$ are 0 and the identity function respectively. We end up with $a \approx -0.9974$, which is reasonably close to the true value $a^\dagger = -1$.

\section{Learning the Poincar\'e normal form of the supercritical Hopf bifurcation}

In this section, we extend our framework to the case of learning  Poincar\'e normal forms of ODEs\footnote{Check Appendix A.4 for more background on Poincar\'e normal forms. }
We extend our Gaussian process (GP) approach to learn a local transformation that maps a trajectory of a nonlinear dynamical system to the corresponding trajectory of its normal form. We will focus on the case of a system with a Hopf bifurcation and plot the true limit cycle of the normal form as well its reconstructed counterpart to show that GPs are effective at learning the map that transforms a system to its Poincar\'e normal form.

\begin{figure}[tbh]
\centering
\subfloat{\label{fig:test9}
\includegraphics[width=.45\textwidth]{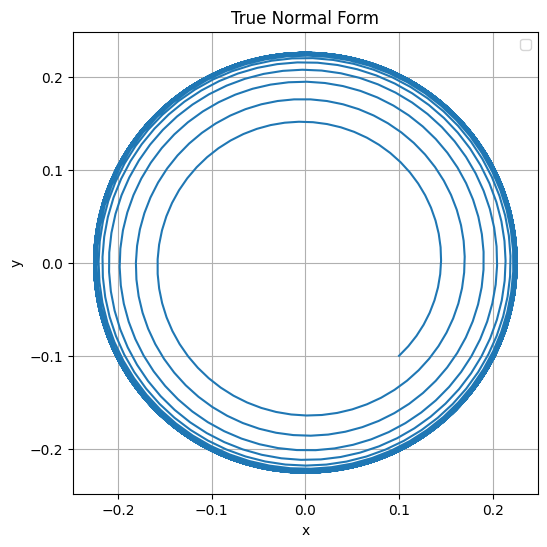}
}
\hfill
\subfloat{\label{fig:test10}
\includegraphics[width=.45\textwidth]{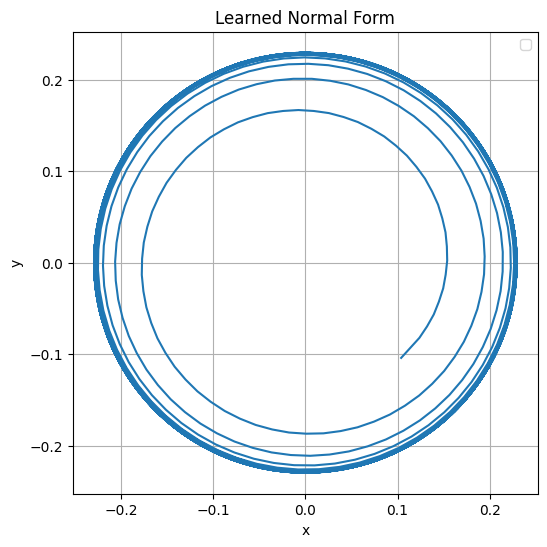}
}
\hfill
\caption{True (left) and learned (right) limit cycles}
\label{limitcycle}
\end{figure}

We examine the case of the Brusselator 

\begin{equation}
\begin{aligned}
\quad & \frac{\text{d}\bar{u}}{\text{d}t} = A+\bar{u}^2\bar{v} -(B+1)\bar{u} \\
\quad & \frac{\text{d}\bar{v}}{\text{d}t} = B\bar{u} -\bar{u}^2\bar{v} \\
\end{aligned},
\end{equation}

where $A$ and $B$ are parameters. The Brusselator undergoes a stable limit cycle for $B > A^2+1$. From the theory of normal forms, there exists a local transformation between the following form of the Brusselator

\begin{equation}
\begin{aligned}
\quad & \frac{\text{d}u}{\text{d}t} = A+(u+A)^2 \left(v+\frac{B}{A} \right) -(B+1)(u+A) \\
\quad & \frac{\text{d}v}{\text{d}t} = B(u+A) -(u+A)^2 \left(v+\frac{B}{A} \right) \\
\end{aligned},
\end{equation}

where we have made the shift of coordinates $u = \bar{u}-A, \; v=\bar{v}-\frac{B}{A}$ to change the equilibrium point to the origin; and the normal form equations, which can be written in Cartesian coordinates:

\begin{equation} \label{nf:hopfcartesian}
\begin{aligned}
    \frac{\text{d}x}{\text{d}t} = (\mu -x^2-y^2)x -y \\
    \frac{\text{d}y}{\text{d}t} = (\mu-x^2-y^2)y + x \\
\end{aligned},
 \end{equation}

 where $\mu$ is a known constant\footnote{From Kuznetsov \cite{kuznetsov}, $\mu = \frac{B-(A^2+1)}{\sqrt{4A^2-(B-A^2-1)^2}}$. In this paper, we run numerical experiments for $A=1.0, B=2.1, \mu \approx 0.05$. For much smaller values of $
 \mu$ around the bifurcation parameter $
 \mu=0$, we expect our method to work. The only drawback when $\mu$ gets smaller is that the radius of the limit cycle, $r \approx \sqrt{\mu}$, gets smaller and the choice of the initial conditions to run our method becomes more sensitive since they have to be within the radius of the limit cycle in order to guarantee that the trajectory coverges to the limit cycle, and the number of collocation points needed gets larger since it takes longer to converge to the limit cycle.},  that is approximately the square of the radius of the limit cycle of our normal form equations; or in polar coordinates (which we shall use in our computations):

\begin{equation} \label{nf:hopf}
    \begin{aligned}
       \quad & \frac{\text{d}r}{\text{d}t} = (\mu-r^2)r\\
       \quad & \frac{\text{d}\theta}{\text{d}t} = 1
    \end{aligned} ,
\end{equation}
 where $x = r \cos \theta$ and $y = r \sin \theta$.

 For a fixed triple $(A,B,\mu)$, the goal is to use CGC to find the GP approximation $H$ of the map $H^\dagger (u,v) = r$. We assume we know that $\theta (t) = t+\theta_0$ because $\theta$ does not depend on $\mu$. Since the solution of the normal form equations is close to the solution of the Brusselator near the origin, $r_0$, the initial condition of the ODE for $r$, is approximately equal to $\sqrt{u^2_0+v^2_0}$ for small values of $u_0$ and $v_0$. The problem of finding the best GP approximation $H$ for a given parameter triple $(A,B,\mu)$ becomes that of finding the $Z$ that minimizes the loss function

\begin{equation}
      \mathcal{L}(H,Z) =  ||H||^2_K + \lambda_1 \mathcal{L}_1 (f,Z) + \lambda_2 \mathcal{L}_2 (Z) + \lambda_3 \big|H(u_0,v_0)-\sqrt{u^2_0+v^2_0} \big|^2,
 \end{equation}

where $Z=(Z_{\cdot,0}, Z_{\cdot,1},Z_{\cdot,2},Z_{\cdot,3},Z_{\cdot,4})$ to be the array of data points approximating the state of the system $(t,u,v,r,\frac{\text{d}r}{\text{d}t})$ at a given iteration step.

Here,

\begin{equation}
    \mathcal{L}_1 (H,Z) = |H (Z_{\cdot,1},Z_{\cdot,2}) - Z_{\cdot,3}|^2 
\end{equation}

is the loss term which forces $Z_{\cdot,3}$, the approximation of $r$, to be the output of $H$ at $(Z_{\cdot,1},Z_{\cdot,2})$ 


\begin{equation}
\mathcal{L}_2 (Z) = |Z_{\cdot,4}-(\mu-Z_{\cdot,3}^2)Z_{\cdot,3}|^2
\end{equation}
 
is the term which enforces the ODE (\ref{nf:hopf}), and

\begin{equation}
 |H (u_0,v_0)-\sqrt{u^2_0+v^2_0}|^2
\end{equation}

is a constraint which maps the initial condition of the Brusselator $(u_0,v_0)$ to $\sqrt{u^2_0+v^2_0}$, which is closed to the true initial condition $r_0$ of the normal form equation for $r$\footnote{We can make such an assumption since the normal form transformation is a near-identity transformation and we can assume that sufficiently small initial conditions of $(u,v)$ will be mapped almost identically to the initial conditions of $(x,y)$, i.e. $(x_0,y_0) \approx (u_0,v_0)$. }. In our simulation, we use $(u_0,v_0) = (0.1,-0.1)$.

 The kernel $K$ is the quartic polynomial kernel $K(\textbf{s},\textbf{t}) = (\textbf{s}^T \textbf{t})^4, \; \textbf{s},\textbf{t} \in \mathbb{R}^2$. We choose a polynomial kernel because both the normal form and the transformations to the normal form are often approximated by Taylor series (polynomial) expansions in the literature. Using the polynomial kernel also allows us to automatically satisfy $H(0,0)=0$ when we use the representer formula for $H$; that is, $H$ maps the origin of the Brusselator to the origin of its normal form. 

\begin{figure}
    \centering
    \includegraphics[width=0.5\linewidth]{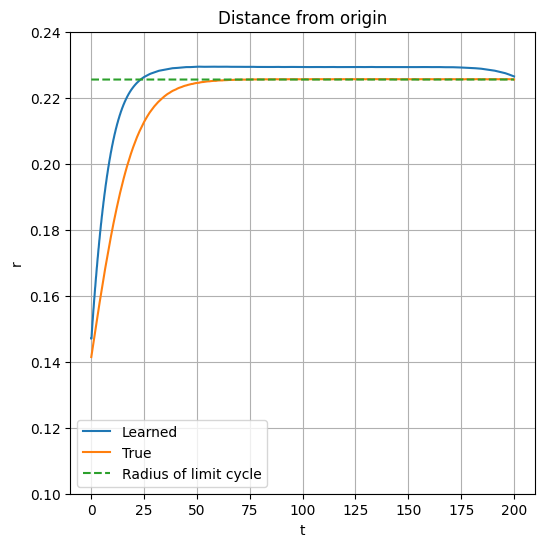}
    \caption{Distance from the origin, $r$, for the true vs. approximate dynamics in polar coordinates}
    \label{fig:rplot}
\end{figure}
\begin{figure}[tbh]
\centering
\subfloat{\label{fig:test11}
\includegraphics[width=.45\textwidth]{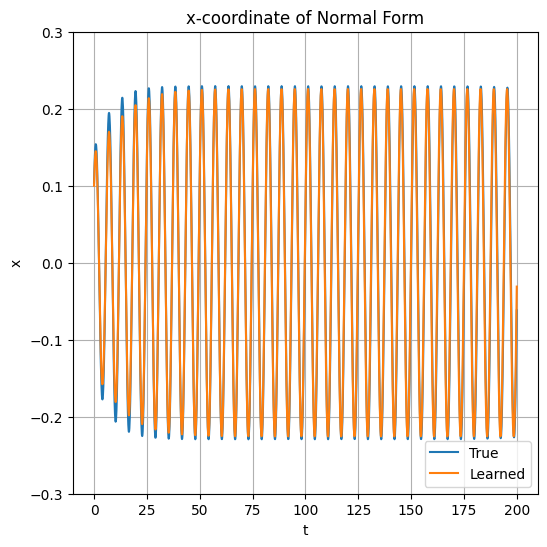}
}
\hfill
\subfloat{\label{fig:test12}
\includegraphics[width=.45\textwidth]{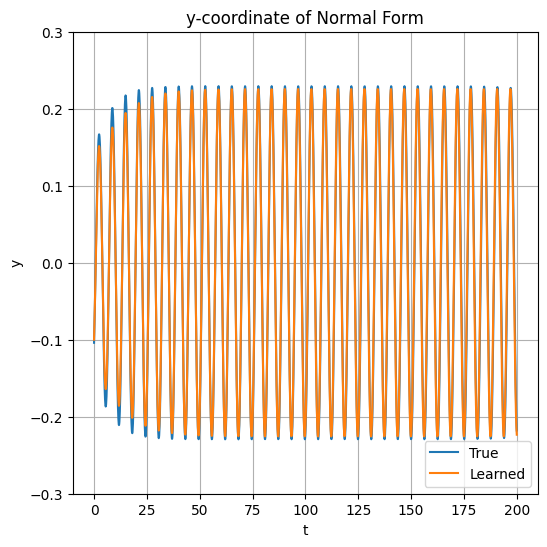}
}
\hfill
\caption{$x$ (left) and $y$ (right) coordinates of normal form equations}
\label{xynormalform}
\end{figure}

We use $N=2000$ points for $(u,v)$ between $t = 0$ and $t=200$ for $(u,v)$ and $A=1$, $B=2.1, \mu = \frac{0.1}{\sqrt{3.99}}$ , $(r_0,\theta_0) = (\frac{\sqrt{2}}{10}, -\frac{\pi}{4})$ in our experiments and we plot true $x = r \cos \theta$ and $y = r \sin \theta$ against our learned $x$ and $y$ in Figure \ref{xynormalform}, as well as the true and learned $r$ as a function of $t$ in Figure \ref{fig:rplot}.

\begin{figure}[tbh]
\centering
\subfloat{\label{fig:test13}
\includegraphics[width=.45\textwidth]{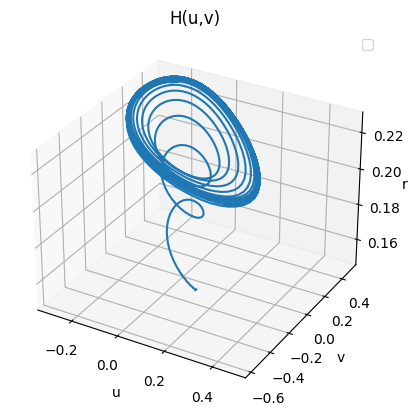}
}
\hfill
\subfloat{\label{fig:test14}
\includegraphics[width=.45\textwidth]{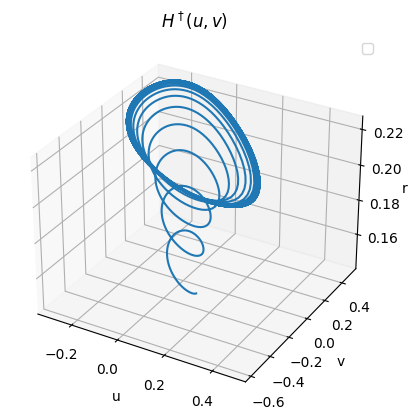}
}
\hfill
\caption{$H$ (left) and $H^\dagger$ (right) as a function of $(u,v)$}
\label{Hnormalform}
\end{figure}

We also display a 3D plot of the true $H^\dagger$ and the learned $H(u,v)$ as a function of $u$ and $v$ in Figure \ref{Hnormalform} to further highlight the similarity between the true and the learned versions of the normal form.

\section{Discussion}

Having learned a couple of known examples of transformations between trajectories of PDEs, we investigate whether this formulation can be applied to other relationships between nonlinear and linear PDEs that have a first-order dependence on time. 

We shall use normal capital letters $P, L$ and $T^\dagger$ to represent maps that take (potentially multiple) derivatives of functions as inputs, and calligraphic letters $\mathcal{P},\mathcal{L}$ and $\mathcal{T}^\dagger$ for operators acting directly on $w$ and $v$.
 
Let $v$ be the solution of a nonlinear PDE which we write in the form 

\begin{equation}
    v_t =  P[(D^\beta u)_{\beta \in J}]:= P[D^{\beta_1} u,...,D^{\beta_n} u] = \mathcal{P}u,
    \label{eq:generalnonlin}
\end{equation}
 where $D^{\beta_n} u ,n\geq 0$ is the $n$th derivative of $u$ with respect to $x$ and $\mathcal{P}$ is a nonlinear operator. For brevity, we shall write $ P[u,Du,...,D^n u]$ as $\mathcal{P}u$. Now let $w$ be a linear PDE which we express as

\begin{equation}
w_t = L[(D^\gamma w)_{\gamma \in \Gamma}]:=L[D w^{\gamma_1},...,D^{\gamma_l} w]  = \mathcal{L} w.
\label{eq:generallin}
\end{equation}

If there exists a transform $\mathcal{T}$ mapping $u_0$ to $w_0$, we write 

\begin{equation}
w_0 =   T^\dagger [(D^{\alpha} u_0)_{\alpha \in I \subset \mathbb{Z}}]:=T^\dagger [D^{\alpha_1} u_0,..., D^{\alpha_m} u_0] = \mathcal{T}^\dagger u_0.
\end{equation}

In other words, the possible inputs of $T^\dagger$ are $u_0$, its derivatives and its antiderivatives.

\begin{figure}[ht]
    \centering
     \includegraphics{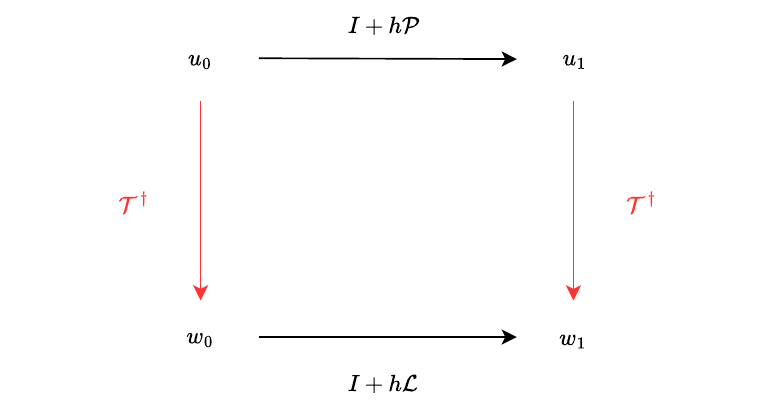}
    \caption{Relationship between nonlinear PDE $u$ and linear PDE $w$}
    \label{GeneralPDEdiagram}
\end{figure}

Given knowledge of the initial conditions $u_0 (x)$ and $w_0(x)$, we can approximate $u_1 (x), w_1 (x)$, the values of $u$ and $w$ at time $t_1 = h$, using the Euler scheme $u_1 = (I + h\mathcal{P})u_0, w_1 = (I+h\mathcal{L})w_0$. Figure \ref{GeneralPDEdiagram} shows that we can reach $w_1$ from $u_0$ in two ways: $w_1 = \mathcal{T}^\dagger u_1 = \mathcal{T}^\dagger \cdot (I + h\mathcal{P})u_0 = (I + h\mathcal{L}) \cdot \mathcal{T}^\dagger u_0$. Taking the limit as $h \rightarrow 0$, we derive the functional differential equation

\begin{equation}
    \frac{\partial \mathcal{T}^\dagger}{\partial \mathcal{P}} \cdot \mathcal{P}= \mathcal{L} \cdot \mathcal{T}^\dagger,
    \label{FDE1}
\end{equation}

where $\frac{\partial \mathcal{T}^\dagger}{\partial \mathcal{P}}$ is the Fréchet derivative of $\mathcal{T}^\dagger$ with respect to $\mathcal{P}$. Then recovering $\mathcal{T}^\dagger$ becomes a question of solving a PDE in terms of $\mathcal{T}^\dagger$.







The optimization problem for $T$, the Gaussian process approximation of $T^\dagger$, now becomes


\begin{equation}
\begin{aligned}
\min_{T \in \mathcal{H}_K} \quad & ||T||^2_K \\
\textrm{s.t.} \quad &     \frac{\partial \mathcal{T}^\dagger}{\partial \mathcal{P}} \cdot \mathcal{P}(x_i)= \mathcal{L} \cdot \mathcal{T}^\dagger (x_i) , \; i = 1, \dots, N\\
\quad & T[D^{\alpha_1} u_0(x_j),..., D^{\alpha_m} u_0 (x_j)] = w(x_j), \; j=1,\dots,l \\
\end{aligned},
\label{eq:colehopfminimizer4}
\end{equation}

where like before, $T[D^{\alpha_1} u_0(x_j),..., D^{\alpha_m} u_0 (x_j)] = w_0 (x_j)$ are arbitrary uniqueness constraints which ensure that the answer to (\ref{eq:colehopfminimizer4}) is not $0$. We remark that in practice, the $u_0 (x_i)$ and $w_0 (x_i)$ may come from multiple pairs of solutions to (\ref{eq:generalnonlin}) and 
 (\ref{eq:generallin}); this will improve the ability of the model to generalize to unseen data and demonstrate that the transform exists in a large region in space.

Then the problem of learning an operator that maps a nonlinear PDE $P$ to a linear PDE $L$ is that of finding the pairs $(P,L)$ such that the RKHS norm of the minimizer $T$ of (\ref{eq:colehopfminimizer4}) remains uniformly bounded independently from the number of constraints $N$ and $l$. As for the number of derivatives and antiderivatives that are required as inputs of the operator $T$, we could start by assuming that $T$ is a function of $u_0$ alone. If there exists no such function $T$ with the property $T(u_0)=w_0$, then the RKHS norm will go to infinity as the number of collocation points increases. In that case, we then investigate whether there exists a function $T(u_0,Du_0) = w_0$ by solving the same minimization problem and repeat the process of increasing the number of inputs. On the other hand, should such a mapping exist, the RKHS norm of $T$ will converge to the norm of $T^\dagger$ as more points are added. We progressively add more inputs as necessary until we find a mapping which is regular.




\section{Code}

The code for this paper can be found at www.github.com/jonghyeon1998/CGC

\section{Conclusion}

We have demonstrated that Gaussian processes (GPs) are capable of recovering the Cole-Hopf transformation from samples of the Burgers equation and its antiderivative, as well as a transformation between the initial conditions of two first-order PDEs.
Although we have mainly focused on an examples of transformations which can be derived exactly, it may be possible to use GPs to learn new transforms which have no closed form using the method outlined in the discussion. The example about normal forms and local transformations between two dynamical systems is also worthy of further investigation. In particular, we shall attempt in future work to learn the radius of the limit cycle of the normal form equations as a function of the parameters of the source dynamical system. 
\appendix
\section{}

\subsection{Reproducing Kernel Hilbert Spaces (RKHS)}

We give a brief overview of reproducing kernel Hilbert spaces as used in statistical learning
theory ~\cite{CuckerandSmale}. Early work developing
the theory of RKHS was undertaken by N. Aronszajn~\cite{aronszajn50reproducing}.

\begin{definition} Let  $\mathcal{H}$  be a Hilbert space of functions on a set ${\mathcal X}$.
Denote by $\langle f, g \rangle$ the inner product on ${\mathcal H}$   and let $\|f\|= \langle f, f \rangle^{1/2}$
be the norm in ${\mathcal H}$, for $f$ and $g \in {\mathcal H}$. We say that ${\mathcal H}$ is a reproducing kernel
Hilbert space (RKHS) if there exists a function $K:{\mathcal X} \times {\mathcal X} \rightarrow \RR$
such that\\
 i. $K_x:=K(X,\cdot)\in {\mathcal{H}}$ for all $x\in {\mathcal{X}}$.\\
ii. $K$ spans ${\mathcal H}$: ${\mathcal H}=\overline{\mbox{span}\{K_x~|~x \in {\mathcal X}\}}$.\\
 iii. $K$ has the {\em reproducing property}:
$\forall f \in {\mathcal H}$, $f(x)=\langle f,K_x \rangle$.\\
$K$ will be called a reproducing kernel of ${\mathcal H}$. ${\mathcal H}_K$  will denote the RKHS ${\mathcal H}$
with reproducing kernel  $K$ where it is convenient to explicitly note this dependence.
\end{definition}

The important properties of reproducing kernels are summarized in the following proposition.
\begin{proposition}\label{prop1} If $K$ is a reproducing kernel of a Hilbert space ${\mathcal H}$, then\\
i. $K(X,y)$ is unique.\\
ii.  $\forall x,y \in {\mathcal X}$, $K(X,y)=K(y,x)$ (symmetry).\\
iii. $\sum_{i,j=1}^q\beta_i\beta_jK(X_i,x_j) \ge 0$ for $\beta_i \in \RR$, $x_i \in {\mathcal X}$ and $q\in\mathbb{N}_+$
(positive definiteness).\\
iv. $\langle K(X,\cdot),K(y,\cdot) \rangle=K(X,y)$.
\end{proposition}
Common examples of reproducing kernels defined on a compact domain $\mathcal{X} \subset \mathrm{R}^n$ are the 
(1) constant kernel: $K(X,y)= k > 0$
(2) linear kernel: $K(X,y)=x\cdot y$
(3) polynomial kernel: $K(X,y)=(1+x\cdot y)^d$ for $d \in \N_+$
(4) Laplace kernel: $K(X,y)=e^{-||x-y||_2/\sigma^2}$, with $\sigma >0$
(5)  Gaussian kernel: $K(X,y)=e^{-||x-y||^2_2/\sigma^2}$, with $\sigma >0$
(6) triangular kernel: $K(X,y)=\max \{0,1-\frac{||x-y||_2^2}{\sigma} \}$, with $\sigma >0$.
(7) locally periodic kernel: $K(X,y)=\sigma^2 e^{-2 \frac{ \sin^2(\pi ||x-y||_2/p)}{\ell^2}}e^{-\frac{||x-y||_2^2}{2 \ell^2}}$, with $\sigma, \ell, p >0$.

\begin{theorem} \label{thm1}
Let $K:{\mathcal X} \times {\mathcal X} \rightarrow \RR$ be a symmetric and positive definite function. Then there
exists a Hilbert space of functions ${\mathcal H}$ defined on ${\mathcal X}$   admitting $K$ as a reproducing Kernel.
Conversely, let  ${\mathcal H}$ be a Hilbert space of functions $f: {\mathcal X} \rightarrow \RR$ satisfying
$\forall x \in {\mathcal X}, \exists \kappa_x>0,$ such that $|f(x)| \le \kappa_x \|f\|_{\mathcal H},
\quad \forall f \in {\mathcal H}. $
Then ${\mathcal H}$ has a reproducing kernel $K$.
\end{theorem}


\begin{theorem}\label{thm4}
 Let $K(X,y)$ be a positive definite kernel on a compact domain or a manifold $X$. Then there exists a Hilbert
space $\mathcal{F}$  and a function $\Phi: X \rightarrow \mathcal{F}$ such that
$$K(X,y)= \langle \Phi(x), \Phi(y) \rangle_{\mathcal{F}} \quad \mbox{for} \quad x,y \in X.$$
 $\Phi$ is called a feature map, and $\mathcal{F}$ a feature space\footnote{The dimension of the feature space can be infinite, for example in the case of the Gaussian kernel.}.
\end{theorem}

\subsection{Function Approximation in RKHSs: An Optimal Recovery Viewpoint} 
In this section, we review function approximation in RKHSs from the point of view of optimal recovery as discussed in \cite{owhadi_scovel_2019}. 

\paragraph{Problem {\bf P}:} Given input/output data $(x_1, y_1),\cdots , (x_N , y_N ) \in \mathcal{X} \times \mathbb{R}$,  recover an unknown function $u^{\ast}$ mapping $\mathcal{X}$ to $\mathbb{R}$ such that
$u^{\ast}(x_i)=y_i$ for $i \in \{1,...,N\}$.

In the setting of optimal recovery, \cite{owhadi_scovel_2019}  Problem {\bf P} can be turned into a well-posed problem by restricting candidates for $u$ to belong to a Banach space of functions $\mathcal{B}$ endowed with a norm $||\cdot||$ and identifying the optimal recovery as the minimizer of the relative error

\begin{equation} \label{game}
    \mbox{min}_v\mbox{max}_u \frac{||u-v||^2}{||u||^2}, 
\end{equation} 
where the max is taken over $u \in \mathcal{B}$ and the min is taken over candidates in $v \in \mathcal{B}$ such that $v(x_i)=u(x_i)=y_i$. For the validity of the constraints $u(x_i) = y_i$,  $\mathcal{B}^{\ast}$, the dual space of $\mathcal{B}$, must contain delta Dirac functions $\phi_i(\cdot)=\delta(\cdot-x_i)$. This problem can be stated as a game between Players I and II and can then be represented as
  \begin{equation}\label{eqdkjdhkjhffORgameban}
\text{\xymatrixcolsep{0pc}\xymatrix{
\text{(Player I)} & u\ar[dr]_{\max}\in \mathcal{B}    &      &v\ar[ld]^{\min}\in L(\Phi,\mathcal{B}) &\text{(Player II)}\\
&&\frac{\|u-v(u)\|}{\|u\|}\,.& &
}}\,
\end{equation}

If $||\cdot||$ is quadratic, i.e. $||u||^2=[Q^{-1}u,u] $ where $[\phi, u]$ stands for the duality product between $\phi \in \mathcal{B}^{\ast}$ and $u \in \mathcal{B}$ and $Q : \mathcal{B}^{\ast}\rightarrow \mathcal{B}$ is a positive symmetric linear bijection (i.e. such that $[\phi, Q \phi] \ge  0$ and $[\psi, Q \phi ] = [\phi, Q \psi]$ for $\phi,\psi \in \mathcal{B}^{\ast} $). In that case the optimal solution of (\ref{game}) has the explicit form 
\begin{equation}\label{sol_rep}v^{\ast}=\sum_{i,j=1}^{N}u(x_i) A_{i,j} Q \phi_j, \end{equation}
where   $A=\Theta^{-1}$ and $\Theta \in \RR^{N \times N}$ is a Gram matrix with entries $\Theta_{i,j}=[\phi_i,Q\phi_j]$.

To recover the classical representer theorem, one defines the reproducing kernel $K$ as $$K(X,y)=[\delta(\cdot-x),Q\delta(\cdot-y)]$$ 
In this case, $(\mathcal{B},||\cdot ||)$ can be seen as an RKHS endowed with the norm
$$||u||^2=\mbox{sup}_{\phi \in \mathcal{B}^\ast}\frac{(\int \phi(x) u(x) dx)^2}{(\int \phi(x) K(X,y) \phi(y) dx dy)}$$
and (\ref{sol_rep}) corresponds to the classical representer theorem 
\begin{equation}\label{eqkjelkjefffhb}
v^{\ast}(\cdot) = y^T AK(X,\cdot),
\end{equation} 
 using the vectorial notation $y^T AK(X,\cdot)=\sum_{i,j=1}^{N}y_iA_{i,j}K(X_j,\cdot)$ with $y_i=u(x_i)$, $A=\Theta^{-1}$ and $\Theta_{i,j} =K(X_i,x_j)$.
  
 Now, let us consider the problem of learning the kernel from data. As introduced in \cite{Owhadi19}, the method of Kernel Flows (KF) is based on the premise that \emph{a kernel is good if there is no significant loss in accuracy in the prediction error if the number of data points is halved}. This led to the introduction of 
 \begin{equation}
     \rho = \frac{||v^{\ast}-v^{s} ||^2}{||v^{\ast} ||^2} 
     \label{eqrhorkhs}
 \end{equation}\footnote{A variant of KF based on Lyapunov exponents, in order to capture long-term behavior of the system, can be found at \cite{BHPhysicaD}, another one based on the Maximum Mean Discrepancy that allows capturing the statistical properties of the system, and that we conjecture to be related to the properties of the Frobenius-Perron operator, can be found at \cite{BHPhysicaD}, and another one based on the Hausdorff distance, that allows reconstructing attractors, can be found at \cite{bh_kfs_p6}. See also \cite{bh_sparse_kfs} for the method of {Sparse Kernel Flows} where we modified (\ref{eqrhorkhs}) by adding an $L_1$ penalty term. An Algorithmic Information Theoretic point of view on KF, and particularly Sparse Kernel Flows,  can be found at \cite{bh_huter_kfs_ait}  where the problem of learning kernels from data is viewed as a problem of data-compression and where it is shown that KF is the ``natural approach'' to learn kernels from data. }
  which is the relative error between 
  $v^\ast$, the optimal recovery \eqref{eqkjelkjefffhb} of $u^\ast$ based on the full dataset
  $X=\{(x_1,y_1),\ldots,(x_N,y_N)\}$, and
  $v^s$  the optimal recovery  of both $u^\ast$ and $v^\ast$ based on half of the dataset $ X^s=\{(x_i,y_i)\mid i \in \mathcal{S}\}$ ($\operatorname{Card}(\mathcal{S})=N/2$) which admits the representation
  \begin{equation}
v^s=(y^s)^T A^s K(X^s,\cdot)
  \end{equation}
 with $y^s=\{y_i\mid i \in \mathcal{S}\}$,
 $x^s=\{x_i\mid i \in \mathcal{S}\}$,
 $A^s=(\Theta^s)^{-1}$, $\Theta^s_{i,j}=K(X_i^s,x_j^s)$.
 This quantity  $\rho$ is directly related to the game in (\ref{eqdkjdhkjhffORgameban}) where one is minimizing the relative error of $v^{\ast}$ versus $v^s$. 
Instead of using the entire the dataset $X$ one may use random subsets $X^{s_1}$ (of $X$) for $v^{\ast}$ and random subsets $ X^{s_2}$ (of $X^{s_1}$) for $v^s$.  In practice, it is computed as   \cite{Owhadi19}

\begin{equation} \label{rho_scalar_s1_s2}
\rho=1-\frac{{Y_{s_2}}^TK(X^{s_2},X^{s_2})^{-1}Y_{s_2}}{{Y_{s_1}}^TK(X^{s_1},X^{s_1})^{-1}Y_{s_1}}
\end{equation}

Writing $\sigma^2(x)=K(X,x)-K(X,X^f)K(X^f,X^f)^{-1}K(X^f,x)$ we have the pointwise error bound
\begin{equation}\label{error_estimate} |u(x)-v^\ast(x)| \leq  \sigma(x) \|u\|_{\Hc},\end{equation}
 Local error estimates such as (\ref{error_estimate}) are
classical in Kriging \cite{Wu92localerror} (see also \cite{owhadi2015bayesian}[Thm. 5.1] for applications to PDEs). $\|u\|_{\Hc}$ is bounded from below (and, in with sufficient data, can be approximated by) by $\sqrt{Y^{f,T} K(X^f,X^f)^{-1} Y^f} $, i.e., the RKHS norm of the interpolant  of $v^\ast$.

\subsection{Derivation of Cole-Hopf}

We outline the theory behind the Cole-Hopf transformation; a more detailed explanation can be found in \cite{evanspde}. Consider the following parabolic nonlinear PDE on $\mathbb{R}$:


\begin{equation}
u_t - \nu u_{xx} + \frac{1}{2} u^2_x = 0
\label{eq:burgers}
\end{equation}

with initial condition $u(x,0)=u_0 (x)$ at time $t=0$ and parameter $\nu \in \mathbb{R}$.

Assume that there exists a smooth function $F : \mathbb{R} \rightarrow \mathbb{R}$ such that $F(u) = w$, where $w(x,t)$ is a solution of a linear PDE.

Then differentiating both sides of $F(u) = w$ with respect to time and space yields the relations


\begin{equation}
 w_t = F'(u) u_t,  \; w_{xx} = F'(u) u_{xx} + F''(u) u^2_x   
 \label{eq:heatode}
\end{equation}

Substituting \ref{eq:heatode} into \ref{eq:burgers} gives us

\begin{equation}
w_t = F'(u) u_t = F'(u) \left( \nu u_{xx} - \frac{1}{2} u^2_x\right) = 
\nu w_{xx}  - \left(\nu F''(u)+\frac{1}{2} F'(u)\right) u^2_x
\label{eq:derivation}
\end{equation}

If we can find a $F$ that solves the ODE $\nu F''+\frac{1}{2} F' = 0$, then \ref{eq:derivation} collapses to

\begin{equation}
w_t = \nu w_{xx}
\end{equation}

, the heat equation. By inspection, the ansatz $F(u) = e^{-\frac{1}{2\nu}u}$ solves the ODE. We can therefore write that

\begin{equation}
w = e^{-\frac{1}{2\nu}u}
\label{eq:antiburgereq}
\end{equation}

where $w$ has initial condition $w(x,0) = e^{-\frac{1}{2\nu}u_0 (x)} := w_0 (x)$. \ref{eq:antiburgereq} is the \emph{Cole-Hopf transformation}.

Since any affine transform $\lambda w+\mu, \; \lambda,\mu\in \mathbb{R}$ of $w$ also solves the heat equation, we can construct a function $D^*$ that maps $u$ to another solution of the heat equation with different initial conditions:

\begin{equation}
    D^* (u):= \lambda e^{-\frac{1}{2\nu}u}+\mu
\end{equation}

Note that we can rearrange \ref{eq:antiburgereq} in terms of $w$ and differentiate with respect to space to obtain a new equation.
\begin{equation}
u_x = v = \mathcal{C}(w,w_x) := -2\nu\frac{w_x}{w} 
\label{eq:colehopf}
\end{equation}

The left hand side $u_x$, which we shall denote hereafter as $v$, is the solution of the viscous Burgers' equation with initial condition $v(x,0)= v_0 (x) $:

\begin{equation}
v_t - \nu v_{xx} + vv_x = 0
\end{equation}

 \subsection{Poincar\'e Normal Forms}
In this section, we briefly review some results on normal forms near equilibria of nonlinear ODEs.

Consider the nonlinear ODE in $\RR^n$ 
\begin{equation}
\label{ode} 
\dot{x}:=F(x)=Ax+f(x),
\end{equation}
with $f \in {C}^{r+1}(\RR^n;\RR^n)$ denoting a vector field on $\RR^n$, $f(0)=0$ and $A=\frac{\partial F}{\partial x}|_{x=0}$. Without loss of generality, we will assume that $A$ is in Jordan form. 

The goal is to find a change of coordinates 
\begin{equation}
\label{chancor} 
x=\xi(y),
\end{equation} with $\xi \in C^r(\RR^n;\RR^n)$ in a neighborhood  of the origin,
such that the Taylor expansion of the vector field after a coordinate transformation (\ref{chancor}) is 
simple in the sense of making essential features of the flow
 of (\ref{ode}) near the equilibrium $x=0$ evident. The resulting vector field, which we will refer to as $g$, is usually called \emph{a normal form}.
 
  The desired simplification of (\ref{ode}) will be
obtained, up to terms of a specified order, by constructing a near identity  coordinate transformation from a
sequence of compositions of coordinate transformations of the form  
(\ref{chancor}) with
\begin{equation}\label{chancor1}
\xi(y)=\exp(\xi^{[k]})(y)=y+\xi^{[k]}(y)+O(|y|^{k+1}),
\end{equation}
 
where $y\in\RR^n$ is close to zero, $\xi^{[k]}\in H_n^k$ ($k \ge 2$), the vector space of
 homogeneous polynomial vector fields of degree $k$ in $n$ variables with values in $\RR^n$, and 
 $\exp(\xi^{[k]})$ denotes the time-one flow of the ODE $\dot{y}=\xi^{[k]}(y)$ for initial conditions in a sufficiently small neighbourhood of the origin $0 \in \RR^n$. 
 We consider a formal power series expansion of $f$ in (\ref{ode}) and write
\begin{equation}
 f(x)=f^{[2]}(x)+f^{[3]}(x)+\ldots,
\end{equation}
with $f^{[k]}\in H_n^k$. From (\ref{chancor1}) we obtain
\begin{equation}\label{chancor2} \xi^{-1}(y)=y-\xi^{[k]}(y)+O(|y|^{2k}). 
\end{equation}
Substituting  (\ref{chancor}), (\ref{chancor1}) and (\ref{chancor2}) in (\ref{ode}), we get
\begin{equation}
\label{y_eqn} \dot y={ A} y+\cdots+f^{[k-1]}(y)
+
f^{[k]}(y)-(L_{A} \xi^{[k]})(y)
+O(|y|^{k+1}),
\end{equation}
with the Lie derivative $L_{A}$ defined on vector fields $f$ as
\begin{equation}
(L_{ A}f)(y):=\frac{\partial f(y)}{\partial y}{A}y-{A}f(y).
\end{equation}
 In the present context $L_A$ is also known as the \emph{homological} operator.

The Lie derivative leaves $H_n^k$ invariant, $L_{\mathcal{A}}:H_n^k\to H_n^k$. We denote its range in $H_n^k$ as 
${\mathcal{R}}^k$ and let ${\mathcal{C}}^k$ denote a complement  of
${\mathcal{R}}^k$ in $H_n^k$ 
\[
\label{eqn:decompesvec}
H_n^k={\mathcal{R}}^k \oplus {\mathcal{C}}^k, \quad k\ge 2. 
\]
We define a \emph{normal form} of $f$ of order $r$ as a Taylor expansion of the vector field with linear part and terms  $g^{[k]}\in{\mathcal{C}}^k$ for $2\leq k\leq r$. 

%
From the preceding it is clear that the elements of the normal form are in ${\mathcal{C}}^k$. The problem of finding the normal form is thus to find the space ${\mathcal{C}}^k$. Essentially, there are two methods  to find the elements of the normal form\footnote{In this paper, we focused on deriving normal forms using the indirect approach. We leave the derivation of normal forms using the PDE formulation to another paper}:
\begin{itemize}
    \item  {\bf Indirect approach ($A$ is semisimple, i.e. diagonalizable over the complex plane)}: doing the computations in ${\mathcal{R}}^k$ by canceling the maximum number of terms and then deducing the elements that do not belong to ${\mathcal{R}}^k$. In this case,  $$H_n^k=\mbox{im}\,L_{\mathcal A}^k \oplus \mbox{ker}\,L_{\mathcal{A}}^k, \quad k\ge 2. $$
\item  {\bf Direct approach}: characterizing the elements in   ${\mathcal{C}}^k$ by solving a PDE. In this case, a different inner product is introduced such that 
$$ H_n^k=\mbox{im}\,L_{\mathcal{A}}^k \oplus \mbox{ker}\,(L_{\mathcal{A}}^k)^{\ast}, \quad k\ge 2. $$

\end{itemize}
One of the advantages of the latter approach is that it extends easily to the case where $A$ is not semisimple.
%
%
%
%
%
\subsubsection{Indirect Approach}
If $A$ is semisimple,  $A$ can be written as $A = \mbox{diag}\{\lambda_1, \cdots,
\lambda_n\}$ where $\lambda_i$, $i=1,\cdots,n$, are the eigenvalues of $A$ then, as in Poincar\'e's original work and its extensions, a normal form of (\ref{ode}) up to order $r \ge 2$ can be chosen so that  its nonlinear part consists of all resonant monomials up to order $r$ where the resonant monomials are \[x_1^{\alpha_1}\cdots x_n^{\alpha_n} e^l\] (with $\sum_{i=1}^n\alpha_i=k \ge 2$, $1 \le j \le n$, and   $ e^{\ell} $ is the $ {\ell}^{th}$
unit vector) are such that the resonance condition
\[\sum_{i=1}^n\lambda_i \alpha_i-\lambda_j=0, \]
is satisfied.
For example,  after quadratic and cubic  changes of
coordinates \begin{eqnarray}
g^{[2]}(y)&=& \sum_{\lambda_i+\lambda_j=\lambda_{\ell}} c_{\ell}^{ij} e^{\ell} y_i y_j\\
g^{[3]}(y)&=& \sum_{\lambda_i+\lambda_j+\lambda_k=\lambda_{\ell}}
c_{\ell}^{ijk} e^{\ell} y_i y_j y_k 
\end{eqnarray}

In this approach, the computations are done on the space ${\mathcal{R}}^k$,  the image of $L_{\mathcal A}^k$, by computing the image of different elements that constitute the basis of $H_n^k$. Those who do not belong to the image of $L_{\mathcal A}^k$ are in its complement and are the resonant terms. 
\subsubsection{Direct Approach} 

A convenient choice of inner product was introduced by Belitskii \cite{belitskii}, Meyer \cite{meyer1} and Elphick \emph{et al.} \cite{elphick}, enabling the characterization of expression of ${\mathcal{C}^k}$ as the kernel of the Lie derivative of $A^\ast$ (the adjoint of linear part $A$ of the vector field at the equilibrium). The normal form
$g$ can be made to satisfy the system of linear PDEs \[L_{A^{\ast}}g^{[k]}=0.\] 
 that can be solved explicitly using the method of characteristics. 
The normal form obtained is some times called am ``inner-product normal form'' \cite{murdock}.  It is based on properly choosing the space complementary to $\mbox{im}(L^k_A)$ such that computations become simple. 
In general, with the above choices made, the normal form is equivariant with respect to the group 
\[G_s=\overline{\{\exp(A^{\ast}_st)~|~t\in\RR\}},\]
i.e. that $x$ is a solution  of the ODE associated to the truncated normal form vector field $\dot{x}=\sum_{k=1}^p g^{[k]}(x)$ if and only if $\gamma(x)$ is also a solution of the same ODE, i.e.
$\frac{d}{dt}(\gamma(x)) =\sum_{k=1}^p g^{[k]}(\gamma(x))$ where $\gamma \in G_s$. The appearance of this equivariance is a simplifying feature. 

\subsection{The Andronov-Hopf Bifurcation }
In this section, we review some results about the  Andronov-Hopf bifurcation. The material is taken from \url{http://www.scholarpedia.org/article/Andronov-Hopf_bifurcation}

\paragraph{Definition}

The Andronov-Hopf bifurcation is a local bifurcation in which a fixed point of a dynamical system loses stability as a pair of complex conjugate eigenvalues cross the imaginary axis of the complex plane. This crossing leads to the appearance or disappearance of a periodic orbit (limit cycle). The bifurcation can be supercritical, where a stable limit cycle is born, or subcritical, where an unstable limit cycle appears.

Consider a dynamical system described by the differential equation:
\[
\dot{x} = f(x,\lambda), \quad x \in \mathbb{R}^n, \quad \lambda \in \mathbb{R},
\]
where \(x\) is the state vector and \(\lambda\) is a parameter. A fixed point \(x_0\) satisfies \(f(x_0, \lambda) = 0\). As the parameter \(\lambda\) varies, the nature of the fixed point may change, leading to different types of bifurcations. The Andronov-Hopf bifurcation occurs when the fixed point undergoes a change in stability due to the movement of a pair of complex conjugate eigenvalues across the imaginary axis.

\paragraph{Two-dimensional Case}

In the case of a two-dimensional system (\(n=2\)), the Andronov-Hopf bifurcation can be described by considering the eigenvalues of the Jacobian matrix at the fixed point. Suppose the system is given by:
\[
\dot{x} = f(x,\lambda), \quad x \in \mathbb{R}^2.
\]
Let \(J(\lambda)\) be the Jacobian matrix evaluated at the fixed point \(x_0(\lambda)\). The eigenvalues \(\mu(\lambda) \pm i\omega(\lambda)\) of \(J(\lambda)\) are complex conjugates. A Hopf bifurcation occurs when: i.) \(\mu(\lambda) = 0\) at \(\lambda = \lambda_0\), ii.) \(\frac{d\mu}{d\lambda} \neq 0\) at \(\lambda = \lambda_0\), iii.) \(\omega(\lambda_0) \neq 0\).

As \(\lambda\) passes through \(\lambda_0\), the real part of the eigenvalues changes sign, leading to a qualitative change in the stability of the fixed point. If \(\mu(\lambda)\) changes from negative to positive, a stable fixed point becomes unstable, and a small stable limit cycle emerges, corresponding to a supercritical Hopf bifurcation. If \(\mu(\lambda)\) changes from positive to negative, the fixed point becomes stable and a small unstable limit cycle disappears, corresponding to a subcritical Hopf bifurcation.

\paragraph{First Lyapunov Coefficient}

Whether the Andronov-Hopf bifurcation is subcritical or supercritical is determined by \(\sigma\), which is the sign of the first Lyapunov coefficient \(l_1(0)\) of the dynamical system near the equilibrium. This coefficient can be computed at \(\alpha = 0\) as follows. Write the Taylor expansion of \(f(x,0)\) at \(x = 0\) as
\[
f(x,0) = A_0x + \frac{1}{2}B(x,x) + \frac{1}{6}C(x,x,x) + O(\|x\|^4),
\]
where \(B(x,y)\) and \(C(x,y,z)\) are the multilinear functions with components
\[
B_j(x,y) = \sum_{k,l=1}^{n} \frac{\partial^2 f_j(\xi,0)}{\partial \xi_k \partial \xi_l} \Bigg|_{\xi=0} x_k y_l,
\]
\[
C_j(x,y,z) = \sum_{k,l,m=1}^{n} \frac{\partial^3 f_j(\xi,0)}{\partial \xi_k \partial \xi_l \partial \xi_m} \Bigg|_{\xi=0} x_k y_l z_m,
\]
where \(j = 1, 2, \dots, n\).

Let \(q \in \mathbb{C}^n\) be a complex eigenvector of \(A_0\) corresponding to the eigenvalue \(i\omega_0\): \(A_0 q = i\omega_0 q\). Introduce also the adjoint eigenvector \(p \in \mathbb{C}^n\): \(A_0^T p = -i\omega_0 p\), \(\langle p, q \rangle = 1\). Here \(\langle p, q \rangle = \overline{p}^T q\) is the inner product in \(\mathbb{C}^n\). Then (see, for example, Kuznetsov (2004))
\[
l_1(0) = \frac{1}{2\omega_0} \text{Re} \left[ \langle p, C(q,q,\overline{q}) \rangle - 2 \langle p, B(q, A_0^{-1}B(q,\overline{q})) \rangle + \langle p, B(\overline{q},(2i\omega_0 I_n - A_0)^{-1}B(q,q)) \rangle \right],
\]
where \(I_n\) is the unit \(n \times n\) matrix. Note that the value (but not the sign) of \(l_1(0)\) depends on the scaling of the eigenvector \(q\). The normalization \(\langle q, q \rangle = 1\) is one of the options to remove this ambiguity. Standard bifurcation software (e.g., MATCONT) computes \(l_1(0)\) automatically.

For planar smooth ODEs with \(x = (u, v)^T\),
\[
f(x,0) = \begin{pmatrix}
0 & \omega_0 \\
-\omega_0 & 0
\end{pmatrix} \begin{pmatrix}
u \\
v
\end{pmatrix} + \begin{pmatrix}
P(u,v) \\
Q(u,v)
\end{pmatrix},
\]
the setting \(q = p = \frac{1}{\sqrt{2}}(1 - i)\) leads to the formula
\[
l_1(0) = \frac{1}{8\omega_0} \left(P_{uuu} + P_{uvv} + Q_{uuv} + Q_{vvv}\right) + \frac{1}{8\omega_0^2} \left[ P_{uv}(P_{uu} + P_{vv}) - Q_{uv}(Q_{uu} + Q_{vv}) - P_{uu}Q_{uu} + P_{vv}Q_{vv} \right],
\]
where the lower indices denote partial derivatives evaluated at \(x = 0\) (cf. Guckenheimer and Holmes, 1983).

The radius \( r \) of the limit cycle that emerges from a supercritical Hopf bifurcation is related to the first Lyapunov coefficient \( l_1(0) \) by the following expression:

\[
r \propto \sqrt{-\frac{l_1(0)}{\omega_0}}
\]

Here, \( \omega_0 \) is the imaginary part of the eigenvalues at the bifurcation point. The sign of \( l_1(0) \) determines whether the bifurcation is supercritical (leading to a stable limit cycle) or subcritical (leading to an unstable limit cycle).

In the case of the Brusselator, $\omega_0 = a > 0$, $\ell_1=-\frac{(2+a^2)}{2a(1+a^2)}$, and 
\[r \propto \frac{\sqrt{2 + a^2}}{a \sqrt{2(1 + a^2)}}\]

\subsection*{Acknowledgments}
BH and YK gratefully acknowledge partial support from the Air Force Office of Scientific Research (award number  FA9550-21-1-0317) and the Department of Energy (award number SA22-0052-S001).
HO gratefully acknowledges partial support from the Air Force Office of Scientific Research under MURI award number FA9550-20-1-0358 (Machine Learning and Physics-Based Modeling and Simulation), from Beyond Limits (Computing Optimal Models), through the JPL Research and Technology Development program award (UQ-aware Machine Learning for Uncertainty Quantification)  and by the Department of Energy under award number DE-SC0023163 (SEA-CROGS: Scalable, Efficient and Accelerated Causal Reasoning Operators, Graphs and Spikes for
Earth and Embedded Systems).  HO and BH acknowledge support from JPL/NASA under the AIST award ``Kernel Flows: Emulating Complex Models for Massive Data Sets''. H.O. is grateful for the Department of Defense Vannevar Bush Faculty Fellowship. 

\bibliographystyle{plain}

\bibliography{missing_dyn, references_p4}

\end{document}